\documentclass[11pt,a4paper]{article}
\usepackage{amsmath, amsthm, amsfonts, amssymb, color}
\usepackage{mathrsfs}
\usepackage{CJKutf8}
\usepackage{color}
\usepackage{bbm}

\newtheorem{thm}{Theorem}[section]

\newtheorem{lem}[thm]{Lemma}
\newtheorem{prp}[thm]{Proposition}

\newtheorem{rem}[thm]{Remark}
\theoremstyle{definition}

\newcommand{\scr}[1]{\mathscr #1}
\definecolor{wco}{rgb}{0.5,0.2,0.3}

\numberwithin{equation}{section} \theoremstyle{remark}

\newcommand{\ua}{\uparrow}

\title{{\bf   Asymptotic log-Harnack inequality for path-distribution dependent SDEs with infinite memory and Dini drift}}
\author{
	{\bf   Xiao-Yu Zhao  }\\
	\footnotesize{ Center for Applied Mathematics and KL-AAGDM, Tianjin University,   300072, China}\\
	\footnotesize{zhxy\_0628@tju.edu.cn}
}
\begin{document}
	\allowdisplaybreaks
	\def\R{\mathbb R}  \def\ff{\frac} \def\ss{\sqrt} \def\B{\mathbf
		B}\def\TO{\mathbb T}
	\def\I{\mathbb I_{\pp M}}\def\p<{\preceq}
	\def\N{\mathbb N} \def\kk{\kappa} \def\m{{\bf m}}
	\def\ee{\varepsilon}\def\ddd{D^*}
	\def\dd{\delta} \def\DD{\Delta} \def\vv{\varepsilon} \def\rr{\rho}
	\def\<{\langle} \def\>{\rangle} \def\GG{\Gamma} \def\gg{\gamma}
	\def\nn{\nabla} \def\pp{\partial} \def\E{\mathbb E}
	\def\d{\text{\rm{d}}} \def\bb{\beta} \def\aa{\alpha} \def\D{\scr D}
	\def\si{\sigma} \def\ess{\text{\rm{ess}}}
	\def\beg{\begin} \def\beq{\begin{equation}}  \def\eed{\end{equation}}\def\F{\scr F}\def\L{\scr L}
	\def\Ric{{\rm Ric}} \def\Hess{\text{\rm{Hess}}}
	\def\e{\text{\rm{e}}} \def\ua{\underline a} \def\OO{\Omega}  \def\oo{\omega}
	\def\tt{\tilde}
	\def\cut{\text{\rm{cut}}} \def\P{\mathbb P} \def\ifn{I_n(f^{\bigotimes n})}
	\def\C{\scr C}      \def\aaa{\mathbf{r}}     \def\r{r}
	\def\gap{\text{\rm{gap}}} \def\prr{\pi_{{\bf m},\varrho}}  \def\r{\mathbf r}
	\def\Z{\mathbb Z} \def\vrr{\varrho} \def\ll{\lambda}
	\def\L{\scr L}\def\Tt{\tt} \def\TT{\tt}\def\II{\mathbb I}
	\def\i{{\rm in}}\def\Sect{{\rm Sect}}  \def\H{\mathbb H}
	\def\M{\scr M}\def\Q{\mathbb Q} \def\texto{\text{o}} \def\LL{\Lambda}
	\def\Rank{{\rm Rank}} \def\B{\scr B} \def\i{{\rm i}} \def\HR{\hat{\R}^d}
	\def\to{\rightarrow}\def\l{\ell}\def\iint{\int}
	\def\EE{\scr E}\def\Cut{{\rm Cut}}\def\W{\mathbb W}
	\def\A{\scr A} \def\Lip{{\rm Lip}}\def\S{\mathbb S}
	\def\BB{\mathbb B}\def\Ent{{\rm Ent}} \def\i{{\rm i}}\def\itparallel{{\it\parallel}}
	\def\g{{\mathbf g}}\def\Sect{{\mathcal Sec}}\def\T{\mathcal T}\def\V{{\bf V}}
	\def\PP{{\bf P}}\def\HL{{\bf L}}\def\Id{{\rm Id}}\def\f{{\bf f}}\def\cut{{\rm cut}}
	\def\Ss{\mathbb S}
	\def\BL{\scr A}\def\Pp{\mathbb P}\def\Pp{\mathbb P} \def\Ee{\mathbb E} \def\8{\infty}\def\1{{\bf 1}}
	\maketitle
	\def\Cc{\mathcal C} \def\t{\theta}
	\begin{abstract} 
		We establish an asymptotic log-Harnack inequality for stochastic differential equations on $\R^d$ whose coefficients depend on the path and distribution over the whole history, allowing the drift to contain a Dini-continuous term. The result is new even in the distribution-independent case.
	\end{abstract}
	\noindent
	Keywords: Path-distribution dependent SDEs,  asymptotic log-Harnack inequality, gradient estimate, Dini continuity.\\
    AMS Subject Classification: 60H10, 60J60, 47G20.
	
	\section{Introduction}
	The dimension-free Harnack inequality with power, first proposed in \cite{W97} for elliptic diffusion semigroups on Riemannian manifolds, has become a powerful analytic tool in probability theory. It has been systematically developed and extended to various stochastic systems, including stochastic differential equations (SDEs), stochastic partial differential equations (SPDEs), path-dependent SDEs, and distribution-dependent SDEs, see \cite{Wbook, RW24} and references therein.
	As a limiting form of the dimension-free Harnack inequality as the power goes to infinity, the log-Harnack inequality was introduced in
	\cite{W10} and has played a key role in the analysis of curvature lower bounds and entropy-cost relations. Further developments have generalized this inequality to metric measure spaces, see \cite{Am, KS}. 
	
	When the noise in a stochastic system is weak, such that the log-Harnack inequality is no longer valid, an effective alternative  is the asymptotic log-Harnack inequality. This has been established,   for instance, in \cite{BWY19} for  path-dependent SDEs and \cite{Xu, HLL20, HLL21} for SPDEs. More recently, this result has been extended to path-distribution dependent SDEs, see \cite{WYZ}.
	This inequality, in particular, yields an asymptotic gradient estimate. 
	In a recent application, \cite{CW25} has employed this inequality in the field of numerical SPDEs to derive numerical weak error bounds for stochastic reaction-diffusion equations near the sharp interface limit.
	
	On the other hand,  path-dependent SDEs with irregular coefficients have also attracted considerable attention.  By combining Zvonkin's transform in \cite{ZY} with truncation arguments, several works have addressed the case where the drift is allowed to be singular only in the present state, since  Zvonkin's transform cannot remove the drift with delay. For instance, see \cite{B18,Ba, H18}.
	
	When singular SDEs depend on both the path and distribution, the analysis becomes substantially more challenging, even for well-posedness, since the distribution is not pathwise determined and  thus cannot be localized through truncation. Under the additional assumption that the path-dependent part of the drift is globally bounded, Huang \cite{H21} established the well-posedness and the Harnack inequality for path-distribution dependent SDEs with finite delay. Without this boundedness condition, well-posedness was proved in \cite{Z24}, and later extended to infinite delay systems in \cite{WYZ}.

	In this paper, we study the  asymptotic log-Harnack inequality   for path-distribution dependent SDEs with infinite memory and Dini drift. The Dini condition is weaker than the Lipschitz assumption and is used below to obtain the regularity of the Zvonkin transformation.
	
	Let $(\R^d, |\cdot |)$ be the $d$-dimensional Euclidean space for some $d\in\mathbb N$.
	Denote by $\R^d\otimes\R^d$   the family of all $d\times d$-matrices
	with real entries, which is  equipped with the operator norm $\|\cdot\|$.    Let  $A^*$ denote the   transpose of $A\in \R^d\otimes \R^d$, and let
	$\|\cdot\|_{\infty}$ be the uniform norm for functions taking values in $\R,\R^d$ or $\R^d\otimes \R^d$.
	Let $\|\cdot\|_{\rm HS}$ denote the Hilbert--Schmidt norm of matrices.
	
	Let $\C:=C((-\8,0];\R^d)$. To describe the path dependence with exponentially decaying memory, for $\tau>0$ set
	\begin{align}\label{CR}
		\C_\tau=\left\{\xi\in\C:\|\xi\|_\tau:=\sup_{s\in(-\infty,0]}(\e^{\tau s}|\xi(s)|)<\8\right\}.
	\end{align}
	This is a complete, but not separable, space under the norm $\|\cdot\|_\tau$. 
	Let $\overline{\C_\tau}$ be its completion under $\|\cdot\|_{\tau+1}$. We equip $\C_\tau$ with the trace Borel $\sigma$-field $\B(\C_\tau)$ inherited from $\overline{\C_\tau}$. Let $\scr P$ and $\scr P(\C_\tau)$ be the sets of all probability measures on $(\R^d,\B(\R^d))$ and $(\C_\tau,\B(\C_\tau))$, respectively, with $\scr P(\C_\tau)$ equipped with the weak topology induced by $\|\cdot\|_{\tau+1}$. Let $\B_b(\C_\tau)$ be the class of bounded measurable functions on $\C_\tau$, and $\B_b^+(\C_\tau)$ the set of strictly positive functions in $\B_b(\C_\tau)$.
	
	Let $(W(t))_{t\geq0}$ be a $d$-dimensional Brownian motion defined on a complete filtered probability space $(\Omega,\F,(\F_t)_{t\geq0},\P)$. Let $X_0=(X(r))_{r\leq0}$ be an $\F_0$-measurable $\C_\tau$-valued random variable.
	We consider the following path-distribution dependent SDE with infinite memory:
	\begin{align}\label{E-1}
		\d X(t)=b(X_t,\L_{X_t})\,\d t+\si(X(t))\,\d W(t), \quad t\ge 0,
	\end{align} 
	where
	for each fixed $t\geq 0$, $X_t(\cdot)\in\C_\tau$ is  defined by
	\begin{equation*}X_t(r):=X(t+r),\ \ r\in(-\infty,0],\end{equation*}
	which is called the segment process of $X(t)$, $\L_{X_t}\in \scr P(\C_\tau)$ is the distribution of $X_t$, and
	$$b:\C_\tau\times\scr P(\C_\tau)\to\R^d,\ \ \
	\si:\R^d\to\R^d\otimes\R^d $$ are measurable. When different probability spaces are involved,
	we use $\L_{X_t|\P}$ in place of $\L_{X_t}$ to emphasize the underlying probability measure.
	
	We can observe that
	\beq\label{EPR} \begin{split}
		&\e^{p\tau t}\|X_t\|_\tau^p:= \sup_{s\in(-\8,0]}\left( \e^{p\tau (t+s)}|X(t+s)|^p\right) \\
		& \leq \sup_{s\in(-\8,0]}\left( \e^{p\tau s}|X(s)|^p\right) +\sup_{s\in[0,t]}\left( \e^{p\tau s}|X(s)|^p\right) \\
		&= \|X_0\|_\tau^p+\sup_{s\in[0,t]}\left( \e^{p\tau s}|X(s)|^p\right),\ \ t\geq0,\,p\geq0.
\end{split}\end{equation}
\subsection{The Wasserstein metric on $\C_\tau$}
{
We focus on the space
$$\scr P_k(\C_\tau):=\big\{\mu\in\scr P(\C_\tau):\|\mu\|_k^k:= \mu(\|\cdot\|_{\tau}^k)<\8\big\},\ \ k>0.$$
For $\mu,\nu\in\scr P_k(\C_\tau)$ and $N\in\mathbb N$, set
$$
\widetilde\W_{k,N}(\mu,\nu):=
\inf_{\pi\in C(\mu,\nu)}
\bigg(\int_{\C_\tau\times\C_\tau}\|\xi-\eta\|_{N,\tau}^k\,\pi(\d\xi,\d\eta)\bigg)^{\ff1{1\vee k}},
$$
where $C(\mu,\nu)$ denotes the set of all couplings of $\mu$ and $\nu$,
$$
\|\xi\|_{N,\tau}:=\sup_{s\in[-N,0]}\e^{\tau s}|\xi(s)|
$$
and  we define
\begin{equation}\label{Wdef}
\W_k(\mu,\nu):=\sup_{N\geq1}\widetilde\W_{k,N}(\mu,\nu).
\end{equation}
Let $\C_{N,\tau}:=C([-N,0];\R^d)$ and let $\mu_N,\nu_N$ be the marginal distributions of $\mu,\nu$ on $\C_{N,\tau}$. For $\rho,\theta\in\scr P_k(\C_{N,\tau})$, set
$$
\W_{k,N}(\rho,\theta):=
\inf_{\pi_N\in C(\rho,\theta)}
\bigg(\int_{\C_{N,\tau}\times\C_{N,\tau}}
\|\xi-\eta\|_{N,\tau}^k\,\pi_N(\d \xi,\d \eta)\bigg)^{\ff1{1\vee k}}.
$$

\begin{prp}\label{PW}
For every $k>0$, $(\scr P_k(\C_\tau),\W_k)$ is a complete metric space, and
\begin{equation}\label{EQ}
\W_{k,N}(\mu_N,\nu_N)=\widetilde\W_{k,N}(\mu,\nu),\qquad N\in\mathbb N.
\end{equation}
Consequently,
\begin{equation}\label{Wsup}
\W_k(\mu,\nu)=\sup_{N\geq1}\W_{k,N}(\mu_N,\nu_N).
\end{equation}
Moreover,
\begin{equation}\label{Wfull}
\W_k(\mu,\nu)=
\inf_{\pi\in C(\mu,\nu)}
\bigg(\int_{\C_\tau\times\C_\tau}
\|\xi-\eta\|_\tau^k\,\pi(\d\xi,\d\eta)\bigg)^{\ff1{1\vee k}},
\end{equation}
and the infimum in \eqref{Wfull} is attained.
\end{prp}

\begin{proof}
Put $q=1\vee k$ and set $c_k=1$ for $0<k\leq1$ and $c_k=2^{k-1}$ for $k>1$. Since $\C_{N,\tau}$ is Polish under $\|\cdot\|_{N,\tau}$, $(\scr P_k(\C_{N,\tau}),\W_{k,N})$ is Polish and an optimal coupling exists. For $0<k<1$, the same conclusion follows by taking $\|\cdot\|_{N,\tau}^k$ as the metric, see \cite[Chapter~6]{Villani}.

We first prove \eqref{EQ}. For any $\pi\in C(\mu,\nu)$, its marginal distribution on $\C_{N,\tau}\times\C_{N,\tau}$ belongs to $C(\mu_N,\nu_N)$. Hence,
\begin{equation*}
\W_{k,N}(\mu_N,\nu_N)\leq\widetilde\W_{k,N}(\mu,\nu).
\end{equation*}
For the converse inequality, let $\pi_N\in C(\mu_N,\nu_N)$ be an optimal coupling and set
\begin{equation*}
r_N:\C_\tau\to\C_{N,\tau},\qquad r_N(\xi)=\xi|_{[-N,0]}.
\end{equation*}
Since each $\|\cdot\|_{M,\tau}$ is continuous with respect to $\|\cdot\|_{\tau+1}$ and
$
\|\cdot\|_\tau=\sup_{M\geq1}\|\cdot\|_{M,\tau},
$
$\C_\tau$ is a Borel subset of $\overline{\C_\tau}$, and hence $(\C_\tau,\B(\C_\tau))$ is a standard Borel space. By the disintegration theorem, see, e.g., \cite{Kallenberg}, $\mu$ and $\nu$ admit regular conditional distributions $K_\mu(x,\d\xi)$ and $K_\nu(y,\d\eta)$ with respect to $r_N$. Thus,
\begin{equation*}
\pi(\d\xi,\d\eta):=
\int_{\C_{N,\tau}\times\C_{N,\tau}}
K_\mu(x,\d\xi)K_\nu(y,\d\eta)\,\pi_N(\d x,\d y)
\end{equation*}
defines a coupling $\pi\in C(\mu,\nu)$. Moreover, $K_\mu(x,\cdot)$ and $K_\nu(y,\cdot)$ are concentrated on $r_N^{-1}(x)$ and $r_N^{-1}(y)$, respectively, for $\mu_N$- and $\nu_N$-a.e. $x,y$. Therefore,
\begin{align*}
\widetilde\W_{k,N}(\mu,\nu)^q
&\leq
\int_{\C_\tau\times\C_\tau}
\|\xi-\eta\|_{N,\tau}^k\,\pi(\d\xi,\d\eta)\\
&=
\int_{\C_{N,\tau}\times\C_{N,\tau}}
\|x-y\|_{N,\tau}^k\,\pi_N(\d x,\d y)
=\W_{k,N}(\mu_N,\nu_N)^q.
\end{align*}
This proves \eqref{EQ}, and hence \eqref{Wsup} holds.

By taking the product coupling and using
\begin{equation*}
\|\xi-\eta\|_{N,\tau}^k\leq c_k\{\|\xi\|_\tau^k+\|\eta\|_\tau^k\},
\end{equation*}
we have $\W_k(\mu,\nu)<\infty$. Moreover, by the triangle inequality for $\W_{k,N}$,
\begin{align*}
\W_k(\mu,\nu)
&=\sup_{N\geq1}\W_{k,N}(\mu_N,\nu_N)\\
&\leq\sup_{N\geq1}\W_{k,N}(\mu_N,\rho_N)
+\sup_{N\geq1}\W_{k,N}(\rho_N,\nu_N)=\W_k(\mu,\rho)+\W_k(\rho,\nu).
\end{align*}
If $\W_k(\mu,\nu)=0$, then $\mu_N=\nu_N$ for every $N$. Since $\C_\tau$ is separable under $\|\cdot\|_{\tau+1}$ and $\B(\C_\tau)$ is generated by the restriction maps $\{r_N\}_{N\geq1}$, this implies $\mu=\nu$. Thus, $\W_k$ is a metric on $\scr P_k(\C_\tau)$.

We next prove \eqref{Wfull}. Since $\|\xi-\eta\|_{N,\tau}\leq\|\xi-\eta\|_\tau$, for every $N\geq1$,
\begin{align*}
\widetilde\W_{k,N}(\mu,\nu)^q
&=\inf_{\pi\in C(\mu,\nu)}
\int_{\C_\tau\times\C_\tau}
\|\xi-\eta\|_{N,\tau}^k\,\pi(\d\xi,\d\eta)\\
&\leq\inf_{\pi\in C(\mu,\nu)}
\int_{\C_\tau\times\C_\tau}
\|\xi-\eta\|_\tau^k\,\pi(\d\xi,\d\eta).
\end{align*}
Taking the supremum over $N$ gives
\begin{equation*}
\W_k(\mu,\nu)^q\leq
\inf_{\pi\in C(\mu,\nu)}
\int_{\C_\tau\times\C_\tau}
\|\xi-\eta\|_\tau^k\,\pi(\d\xi,\d\eta).
\end{equation*}

For the converse inequality, for every $N\geq1$ take $\pi_N\in C(\mu,\nu)$ such that
\begin{equation*}
\int_{\C_\tau\times\C_\tau}
\|\xi-\eta\|_{N,\tau}^k\,\pi_N(\d\xi,\d\eta)
\leq\widetilde\W_{k,N}(\mu,\nu)^q+N^{-1}
\leq\W_k(\mu,\nu)^q+N^{-1}.
\end{equation*}
Regarding $\mu,\nu$ and $\pi_N$ as measures on the corresponding Polish completions, the fixed marginals imply that $\{\pi_N\}_{N\geq1}$ is tight. Hence, by Prokhorov's theorem, see, e.g., \cite{Kallenberg}, there exist a subsequence $\{N_j\}_{j\geq1}$ with $\lim_{j\to\infty}N_j=\8$ and a probability measure $\pi$ on $\overline{\C_\tau}\times\overline{\C_\tau}$ such that
\begin{equation*}
\lim_{j\to\infty}
\int F\,\d\pi_{N_j}
=
\int F\,\d\pi,
\qquad
F\in C_b(\overline{\C_\tau}\times\overline{\C_\tau}).
\end{equation*}
Since the coordinate projections are continuous, $\pi$ has marginals $\mu$ and $\nu$. Moreover, $\mu(\C_\tau)=\nu(\C_\tau)=1$, so that $\pi(\C_\tau\times\C_\tau)=1$. Thus, $\pi\in C(\mu,\nu)$.

For any fixed $M\geq1$, the map $(\xi,\eta)\mapsto\|\xi-\eta\|_{M,\tau}^k$ extends to a non-negative continuous function on $\overline{\C_\tau}\times\overline{\C_\tau}$. Hence, letting $N_j\geq M$ for all sufficiently large $j$,
\begin{align*}
&\int_{\C_\tau\times\C_\tau}
\|\xi-\eta\|_{M,\tau}^k\,\pi(\d\xi,\d\eta)\leq\liminf_{j\to\infty}
\int_{\C_\tau\times\C_\tau}
\|\xi-\eta\|_{M,\tau}^k\,\pi_{N_j}(\d\xi,\d\eta)\\
&\quad\leq\liminf_{j\to\infty}
\int_{\C_\tau\times\C_\tau}
\|\xi-\eta\|_{N_j,\tau}^k\,\pi_{N_j}(\d\xi,\d\eta)\leq\W_k(\mu,\nu)^q.
\end{align*}
Since $\|\xi-\eta\|_{M,\tau}^k$ is increasing in $M$ and
$
\lim_{M\to\infty}\|\xi-\eta\|_{M,\tau}^k
=
\|\xi-\eta\|_\tau^k,
$
the monotone convergence theorem yields
\begin{equation*}
\int_{\C_\tau\times\C_\tau}
\|\xi-\eta\|_\tau^k\,\pi(\d\xi,\d\eta)
\leq\W_k(\mu,\nu)^q.
\end{equation*}
Together with the reverse inequality, this proves \eqref{Wfull}. Moreover,
\begin{equation*}
\int_{\C_\tau\times\C_\tau}
\|\xi-\eta\|_\tau^k\,\pi(\d\xi,\d\eta)
=\W_k(\mu,\nu)^q,
\end{equation*}
so the infimum in \eqref{Wfull} is attained.

It remains to prove the completeness. Let $\{\mu^{(n)}\}_{n\geq1}$ be a $\W_k$-Cauchy sequence. 
\begin{equation*}
\W_{k,N}(\mu_N^{(n)},\mu_N^{(m)})\leq\W_k(\mu^{(n)},\mu^{(m)}),\ \ N\geq1,
\end{equation*}
so that $\{\mu_N^{(n)}\}_{n\geq1}$ is $\W_{k,N}$-Cauchy. Hence, $\mu_N^{(n)}$ has a unique limit $\mu_N$ in
$\scr P_k(\C_{N,\tau})$ under $\W_{k,N}$, and
$\{\mu_N\}_{N\geq1}$ is consistent. By Kolmogorov's extension theorem, see, e.g., \cite{Kallenberg}, there exists a probability measure $\bar\mu$ on $C((-\infty,0];\R^d)$ with $\{\mu_N\}_{N\geq1}$ as marginal distributions.

We next show that $\bar\mu$ is concentrated on $\C_\tau$. Choose
$n_0\geq1$ such that
\begin{equation*}
	\W_k(\mu^{(n)},\mu^{(n_0)})\leq1,\qquad n\geq n_0.
\end{equation*}
Using an optimal coupling of $\mu_N^{(n)}$ and $\mu_N^{(n_0)}$ and
the inequality
\begin{equation*}
	\|x\|_{N,\tau}^k
	\leq c_k\big\{\|x-y\|_{N,\tau}^k+\|y\|_{N,\tau}^k\big\},
\end{equation*}
 we obtain
\begin{equation*}
	\sup_{N\geq1,\,n\geq n_0}
	\mu_N^{(n)}(\|\cdot\|_{N,\tau}^k)
	\leq
	c_k\left\{1+\mu^{(n_0)}(\|\cdot\|_\tau^k)\right\}.
\end{equation*}
Since $\mu_N^{(n)}$ converges weakly to $\mu_N$, the Portmanteau theorem yields
\begin{equation*}
	\sup_{N\geq1}\mu_N(\|\cdot\|_{N,\tau}^k)
	\leq
	c_k\left\{1+\mu^{(n_0)}(\|\cdot\|_\tau^k)\right\}<\infty.
\end{equation*}
Since $\mu_N$ is the $[-N,0]$-marginal of $\bar\mu$ and
$
	\lim_{N\to\infty}\|\xi\|_{N,\tau}^k=\|\xi\|_\tau^k,
$
the monotone convergence theorem gives
\begin{equation*}
	\bar\mu(\|\cdot\|_\tau^k)
	=
	\lim_{N\to\infty}
	\mu_N(\|\cdot\|_{N,\tau}^k)
	<\infty.
\end{equation*}
Thus, $\bar\mu(\C_\tau)=1$, and hence $\bar\mu$ defines a measure
$\mu\in\scr P_k(\C_\tau)$.

Finally,
\begin{align*}
\W_k(\mu^{(n)},\mu)
&=\sup_{N\geq1}\W_{k,N}(\mu_N^{(n)},\mu_N)=\sup_{N\geq1}\lim_{m\to\infty}\W_{k,N}(\mu_N^{(n)},\mu_N^{(m)})\\
&\leq\liminf_{m\to\infty}\sup_{N\geq1}\W_{k,N}(\mu_N^{(n)},\mu_N^{(m)})=\liminf_{m\to\infty}\W_k(\mu^{(n)},\mu^{(m)}).
\end{align*}
Since $\{\mu^{(n)}\}_{n\geq1}$ is $\W_k$-Cauchy,
\begin{equation*}
\lim_{n\to\infty}\W_k(\mu^{(n)},\mu)=0.
\end{equation*}
This proves the completeness.
\end{proof}
}
\section{Main Results}

To establish the asymptotic log-Harnack inequality, we decompose $b$ as
$$b(\xi,\mu)= b^{(0)}(\xi(0))+ b^{(1)}(\xi,\mu),\ \ \xi\in \C_\tau,\ \mu\in \scr P_2(\C_\tau),$$ 
{where $b^{(0)}$ may be Dini continuous, namely uniformly continuous with a modulus of continuity belonging to the class}
\begin{align}\label{D}
\D:=\left\{\varphi:[0,\8)\to[0,\8) \text{ is increasing and concave, } \varphi(0)=0, \int_0^1\ff{\varphi(s)}s\,\d s<\8\right\}.
\end{align} 
{The Dini continuity of $b^{(0)}$ is used to ensure the regularity of the transformed diffusion coefficient in the Zvonkin transformation. We explain in the Remark \ref{dini} (1) why this property is needed for the present coupling argument.}

 For any $\xi\in\C_\tau$, let
$\xi^0\in \C_\tau$ be defined as
$$\xi^0(r)=\xi(0),\ \ r\le 0.$$ 
{We impose the following assumptions on the diffusion coefficient, the path-distribution dependent part of the drift, and the Dini-continuous part $b^{(0)}$.}

\begin{enumerate}
\item[($H$)]There exist constants $K,K_1>0$, $\alpha\in[0,1]$ and a function $\varphi\in\D$ such that $b=b^{(0)}+b^{(1)}$  and $a:=\si\si^*$ satisfy the following conditions.
\begin{enumerate}
\item[($H_1$)]   $a:=\si\si^*$ is invertible with $\| a\|_\8+\| (a)^{-1}\|_\8\leq K$.
\item[($H_2$)] For any $\xi,\eta\in\C_\tau$ and $\mu,\nu\in\scr P_2(\C_\tau)$,
\begin{align}\label{blip}
&  |b^{(1)}(\xi,\mu)-b^{(1)}(\eta,\nu)|\leq K\|\xi-\eta\|_\tau+K_1\W_{2}(\mu,\nu),\\
&  |b^{(1)}(\xi,\mu){-b^{(1)}(\xi^0,\mu)}|\leq {K}\left\{1+\|\xi\|_\tau^\alpha\right\}+K_1\|\mu\|_2\notag.
\end{align}
\item[$(H_3)$] $\|b^{(0)}\|_{\8}<\8$, and for any $x,y\in\R^d$,  
\begin{equation*}
\begin{split}
	|b^{(0)}(x)-b^{(0)}(y)|\leq\varphi(|x-y|),\ \ \|\si(x)-\si(y)\|_{HS}\leq K|x-y|.
\end{split}
\end{equation*}
\end{enumerate}
\end{enumerate}

In accordance with \cite[Theorem 2.3]{WYZ} and Young's inequality, the above 
assumptions guarantee that \eqref{E-1} is well-posed for distributions in $\scr P_{k,\e}^\alpha(\C_\tau),$ $k\geq2$, where 
\begin{align*}
\scr P_{k,\e}^\alpha(\C_\tau):=\left\{\mu\in\scr P_k(\C_\tau): \mu(\e^{\dd\|\cdot\|_\tau^{2\alpha}})<\infty \text{ for some } \dd\in(0,1)\right\}.
\end{align*}
In the distribution-dependent estimates below, the initial laws are fixed, and generic constants may depend, besides the structural parameters, on their finite exponential-moment bounds.
In this case, for any  $\xi\in\C_{\tau}$ with $\L_{\xi}\in\scr P_{k,\e}^\alpha(\C_\tau),$ let $P_t^*\gamma=\L_{X_t^\xi}$ with $\gamma=\L_{\xi}$, and
\begin{align*}
P_tf(\gamma):=\E[f(X_t^\xi)]=\int_{\C_\tau}f\,\d\left\{P_t^*\gamma\right\},\ \ t\geq0,\,f\in\B_b(\C_\tau).
\end{align*}
We first consider the distribution-independent case of \eqref{E-1} with Dini drift, that is, $K_1=0$ in $(H_2)$. In this case, for $\xi\in\C_\tau$, write
\begin{equation*}
P_tf(\xi):=\E[f(X_t^\xi)],\qquad f\in\B_b(\C_\tau),
\end{equation*}
and call $\mu$ an invariant probability measure of $P_t$ if
\begin{equation*}
\int_{\C_\tau}P_tf(\xi)\,\mu(\d\xi)=\int_{\C_\tau}f(\xi)\,\mu(\d\xi),
\qquad f\in\B_b(\C_\tau).
\end{equation*}
We establish the following asymptotic log-Harnack inequality.
\begin{thm}\label{T01}
Assume $(H)$ with $K_1=0$. Then, for any $\tau_0\in(0,\tau)$, there exists a constant $c>0$ such that
\begin{align}\label{alog2}
P_t\log f(\eta)
\leq&\ \log P_tf(\xi)
+c(1+\|\eta\|_\tau^{2\alpha})\|\xi-\eta\|_\tau^2\notag\\
&+c\e^{-\tau_0t}\|\nn\log f\|_\8\|\xi-\eta\|_\tau
\end{align}
holds for any $t\geq0$, $\xi,\eta\in\C_\tau$ and $f\in\B_b^+(\C_\tau)$ with $\|\nn\log f\|_\8<\8$.
\end{thm}
Below,
we write $f\in C_{b,L}(\C_\tau)$ if $$\|f\|_\infty+\|\nn f\|_\infty<\infty,$$
where $$|\nn f|(\xi):=\limsup_{\eta\to\xi}\ff{|f(\xi)-f(\eta)|}{\|\xi-\eta\|_\tau}$$
is the Lipschitz constant of  $f$ at point $\xi\in\C_\tau$, which coincides with the   norm of the gradient   $\nn f(\xi)$ if $f$ is G\^ateaux differentiable at $\xi$.

The next result is an immediate consequence of Theorem \ref{T01} and \cite[Theorem 2.1]{BWY19}. For $\xi,\eta\in\C_\tau$, take $E=\C_\tau$, $\rho(\xi,\eta)=\|\xi-\eta\|_\tau$, and
\begin{align*}
\Gamma_t(\xi,\eta)&=c\e^{-\tau_0t},\\
\LL(\xi,\eta)&=c\{2+\|\xi\|_\tau^{2\alpha}+\|\eta\|_\tau^{2\alpha}\}.
\end{align*}
Write $\Gamma_t(\xi)=\Gamma_t(\xi,\xi)$ and $\LL(\xi)=\LL(\xi,\xi)$.
\beg{cor} In the situation of Theorem $\ref{T01}$,
the following assertions hold.
\beg{enumerate}
\item[$(1)$] For any $t\geq0,\xi\in\C_\tau$ and $f\in C_{b,L}(\C_\tau)$,
\begin{equation*}
|\nn P_t f|(\xi)
\leq
\ss{2\LL(\xi)}
\ss{[P_tf^2(\xi)-(P_tf(\xi))^2]}
+\|\nn f\|_\infty\Gamma_t(\xi).
\end{equation*}
Moreover, $P_t$ is asymptotically strong Feller.
\item[$(2)$] If {$P_t$ has an invariant probability measure $\mu$,} then
\begin{align*}
\limsup_{t\to\8}P_t f(\xi)
\leq
\log\left(
\ff{\mu(\e^f)}
{\int_{\C_\tau}\e^{-\LL(\xi,\eta)\|\xi-\eta\|_\tau^2}\mu(\d\eta)}
\right),
\quad
\xi\in\C_\tau,\quad
f\in\B_b^+(\C_\tau)\cap C_{b,L}(\C_\tau).
\end{align*}
Consequently, for any closed set $A\subset\C_\tau$ with $\mu(A)=0$,
\begin{equation*}
\lim_{t\to\8}P_t\1_A(\xi)=0,\qquad \xi\in\C_\tau.
\end{equation*}
\item[$(3)$] Let $\xi\in\C_\tau$ and $A\subset\C_\tau$ be a measurable set such that
\begin{equation*}
\delta(\xi,A):=\liminf_{t\to\8}P_t(\xi,A)>0.
\end{equation*}
{For any $\vv>0$, let
\begin{equation*}
A_\vv:=\{\eta\in\C_\tau:\inf_{\zeta\in A}\|\eta-\zeta\|_\tau<\vv\}.
\end{equation*}}
Then
\begin{equation*}
\liminf_{t\to\8}P_t(\eta,A_\vv)>0,\qquad \eta\in\C_\tau.
\end{equation*}
More precisely, for any $\vv_0\in(0,\delta(\xi,A))$, there exists a constant $t_0>0$ such that for any $\eta\in\C_\tau$ and $\vv>0$,
\begin{align*}
\inf\bigg\{
P_t(\eta,A_\vv):
t>t_0\vee
\Big(
\ff1{\tau_0}
\log
\ff{\Gamma_0(\xi,\eta)\|\xi-\eta\|_\tau}{\vv\vv_0}
\Big)
\bigg\}>0.
\end{align*}
\end{enumerate}
\end{cor}
Building on Theorem \ref{T01}, we extend the asymptotic log-Harnack inequality to the path-distribution dependent \eqref{E-1} with {infinite memory} and Dini drift.
\begin{thm}\label{T02}
{Assume $(H)$. Then the following assertions hold:}
\begin{enumerate}
\item[(1)] {For any $\vv(\alpha)$ with $\vv(0)=0$ and $\vv(\alpha)>0$ for $\alpha>0$, and any $\mu,\nu\in\scr P_{2+\vv(\alpha),\e}^\alpha(\C_\tau)$, there exists an increasing locally bounded function $c_0:[0,\8)\to(0,\8)$ such that}
\begin{align}\label{wt}
\W_2(P_t^*\mu,P_t^*\nu)
\leq
c_0(t)\W_{2+\vv(\alpha)}(\mu,\nu),\qquad t\geq0.
\end{align}
\item[(2)] For any $\tau_0\in(0,\tau)$, any $\vv(\alpha)$ as in $(1)$, and any $\mu,\nu\in\scr P_{2+\vv(\alpha),\e}^\alpha(\C_\tau)$, {there exists a constant $c>0$ such that}
\begin{equation}\label{log1}
\begin{aligned}
P_t\log f(\nu)\leq&
\log P_tf(\mu)+c_1(t)\W_{2+\vv(\alpha)}(\mu,\nu)^2+c\e^{-\tau_0t}\|\nn\log f\|_\8\W_{2+\vv(\alpha)}(\mu,\nu),
\end{aligned}
\end{equation}
holds for any $t\geq0$ and {$f\in\B_b^+(\C_\tau)$} with $\|\nn\log f\|_\8<\8$, where
\begin{equation}\label{c1t}
{
c_1(t):=c\left\{1+K_1^2\int_0^t c_0(s)^2\,\d s\right\}
\leq c\{1+K_1^2t\,c_0(t)^2\},\qquad t\geq0.}
\end{equation}
Consequently, for any $t\geq0$ and $f\in C_{b,L}(\C_\tau)$,
\begin{align}\label{lip}
|P_tf(\mu)-P_tf(\nu)|
\leq
\W_{2+\vv(\alpha)}(\mu,\nu)
\left[
c\e^{-\tau_0t}\|\nn f\|_\8
+2\ss{c_1(t)}\|f\|_\8
\right].
\end{align}
\end{enumerate}
\end{thm}

\begin{rem}
The term involving $K_1$ in \eqref{c1t} comes from the first change of measure in the proof of \eqref{log1}. When $K_1=0$, this change of measure disappears. Taking $\mu=\dd_\xi$ and $\nu=\dd_\eta$ and keeping the dependence on the initial segment as in the proof of Theorem \ref{T01} gives \eqref{alog2}.
\end{rem}
%
%

\section{Proof of Theorem \ref{T01}}
In this section, we consider the following path-distribution independent SDE with infinite memory on $\R^d$:
\beq\label{E-2}
\d X(t)= b(X_t)\mathrm{d} t +\si(X(t))\mathrm{d} W(t),\ \  t>0,\ \ X_0=\xi\in\C_\tau,
\end{equation}
and the associated Markov semigroup
$$	P_tf(\xi):=\E \big[f(X_t^\xi)\big], \quad\ \  t \ge 0,\,\xi\in \C_\tau,\,f\in \B_b(\C_\tau).
$$  
\paragraph{$(a_1)$}We first use Zvonkin's transformation to deal with the non-Lipschitz drift in \eqref{E-2}. {For $x\in\R^d$, let $x^0\in\C_\tau$ be the constant path $x^0(r)=x$, $r\leq0$, and define $\bar b:\R^d\to\R^d$ and $\tt b:\C_\tau\to\R^d$ by
\begin{equation*}
\bar b(x):=b^{(1)}(x^0),\qquad
\tt b(\xi):=b^{(1)}(\xi)-\bar b(\xi(0)).
\end{equation*}
Then
\begin{align}\label{btilde}
|\bar b(x)-\bar b(y)|\leq K|x-y|, \ \
|\tt b(\xi)-\tt b(\eta)|\leq2K\|\xi-\eta\|_\tau,\ \
|\tt b(\xi)|\leq K(1+\|\xi\|_\tau^\alpha).
\end{align}}
Let $u:\R^d\to\R^d$ solve the elliptic equation
\begin{align}\label{ee1}
\left(
\nn_{b^{(0)}+\bar b}
+\ff12{\rm tr}\{a\nn^2\}
-\lambda
\right)u=-b^{(0)}.
\end{align}
{Here $\nn_vu:=\langle v,\nn u\rangle$.}
Due to assumption $(H)$, the same Zvonkin regularity estimate as in \cite[(2.18)]{RW23}, with the Lipschitz first-order term treated as in \cite{ZY}, gives a constant $c_1>0$ such that
\begin{align}\label{nu1}
\|\nn^2 u\|_\8\leq c_1,
\end{align}
and {we choose $\lambda>0$ sufficiently large such that}
\begin{align}\label{nu2}
\|u\|_\8+\|\nn u\|_\8\leq\ff12.
\end{align}
Define
\begin{equation}
\begin{aligned}\label{Theta}
&\Theta(x)=x+u(x), \ \ x\in\R^d,\\
&(\Theta^\#(\xi))(\t)=\Theta(\xi(\t)),\ \ \t\in(-\8,0],\,\xi\in\C_\tau.
\end{aligned}
\end{equation}
Obviously, $\Theta^\#:\C_\tau\to\C_\tau$ is invertible with
\begin{equation*}
\{(\Theta^\#)^{-1}(\xi)\}(\t)
=\Theta^{-1}(\xi(\t)),
\quad
\xi\in\C_\tau,\ \t\in(-\8,0].
\end{equation*}
Moreover, for a strong solution to \eqref{E-2} with $X_0=\xi\in\C_\tau$, using It\^o's formula and \eqref{ee1}, $Y^\xi(t)=\Theta(X^\xi(t))$ solves
\begin{align}\label{dY}
\d Y^\xi(t)=\hat b(Y^\xi_t)\,\d t+\hat\si(Y^\xi(t))\,\d W(t),\ \ t\geq0,
\end{align}
with $Y^\xi_0=\Theta^\#(\xi)$, where
\begin{equation}
\begin{split}\label{tt}
&
\hat b(\xi)
=
\bar b(\Theta^{-1}(\xi(0)))+\lambda u(\Theta^{-1}(\xi(0)))
+\nn\Theta(\Theta^{-1}(\xi(0)))\tt b((\Theta^\#)^{-1}(\xi)),\\
&
\hat\si(x)
=
\{\nn\Theta(\cdot)\si(\cdot)\}(\Theta^{-1}(x)).
\end{split}
\end{equation}
Thus, $\hat X^\xi(t):=Y^{(\Theta^\#)^{-1}(\xi)}(t)$ solves \eqref{dY} with $\hat X^\xi_0=\xi$.
By $(H)$ with $K_1=0$, \eqref{btilde}, \eqref{nu1} and \eqref{nu2}, there exists a constant $c_2>0$ such that
\begin{equation}\label{ttl}
\begin{split}
&|\hat b(\xi)-\hat b(\eta)|
\leq c_2\|\xi-\eta\|_\tau
+c_2\|\eta\|_\tau^\alpha|\xi(0)-\eta(0)|,
\qquad \xi,\eta\in\C_\tau,\\
&\|\hat\si(x)-\hat\si(y)\|_{HS}
\leq c_2|x-y|,
\qquad x,y\in\R^d, \ \
\|\hat\si\|_\8+\|\hat\si^{-1}\|_\8\leq c_2.
\end{split}
\end{equation}
Then the proof of \cite[Theorem 2.1]{WYZ} yields that SDE \eqref{dY} has a unique {non-explosive} strong solution for any initial value $\xi\in\C_\tau$. Let $\hat P_t$ be the semigroup associated with $\hat X^\xi(t)$. We observe that it suffices to consider the asymptotic log-Harnack inequality for $\hat P_t$.

Indeed, since
\begin{equation*}
X^\xi_t
=(\Theta^\#)^{-1}(Y^\xi_t)
=(\Theta^\#)^{-1}\left(\hat X^{\Theta^\#(\xi)}_t\right),
\end{equation*}
we have
\begin{align*}
P_tf(\xi)
&=\E f(X_t^\xi)
=\hat P_t\left(f\circ(\Theta^\#)^{-1}\right)(\Theta^\#(\xi)).
\end{align*}
Once the asymptotic log-Harnack inequality for $\hat P_t$ is established below, \eqref{nu2} implies that there exist constants $c_3,c_4>0$ such that, for any $t\geq0$, $\xi,\eta\in\C_\tau$ and $f\in\B_b^+(\C_\tau)$ with $\|\nn\log f\|_\8<\8$,
\begin{equation}\label{alp}
\begin{aligned}
&P_t\log f(\eta)-\log P_t f(\xi)\\
&={
\hat P_t\left(\log\left(f\circ(\Theta^\#)^{-1}\right)\right)(\Theta^\#(\eta))
-\log\hat P_t\left(f\circ(\Theta^\#)^{-1}\right)(\Theta^\#(\xi))}\\
&\leq
c_3\{1+\|\Theta^\#(\eta)\|_\tau^{2\alpha}\}
\|\Theta^\#(\xi)-\Theta^\#(\eta)\|_\tau^2\\
&\quad+c_3\e^{-\tau_0t}\|\nn\log f\|_\8
\|\Theta^\#(\xi)-\Theta^\#(\eta)\|_\tau\\
&\leq
c_4(1+\|\eta\|_\tau^{2\alpha})\|\xi-\eta\|_\tau^2
+c_4\e^{-\tau_0t}\|\nn\log f\|_\8\|\xi-\eta\|_\tau.
\end{aligned}
\end{equation}

\paragraph{$(a_2)$}Below, we only consider the regular SDE \eqref{dY}.
For any $\xi\in\C_\tau$, simply denote $\hat X_t=\hat X^\xi_t,\,\hat X(t)=\hat X^\xi(t)$.
Let $\kappa>0$ be chosen sufficiently large below and consider
\begin{equation}\label{ttY}
\begin{aligned}
\d\hat Y(t)=&
{\Big\{
\hat b(\hat Y_t)
+\hat\si(\hat Y(t))\hat\si(\hat X(t))^{-1}
\big[
\hat b(\hat X_t)-\hat b(\hat Y_t)
+\kappa(\hat X(t)-\hat Y(t))
\big]
\Big\}}\,\d t\\
&+\hat\si(\hat Y(t))\,\d W(t),
\qquad t\geq0,\quad \hat Y_0=\eta\in\C_\tau.
\end{aligned}
\end{equation}
By assumption $(H)$, \eqref{ttl} and the proof of \cite[Theorem 2.1]{WYZ}, \eqref{ttY} has a unique strong segment solution $(\hat Y_t)_{t\geq0}$.
Let
\begin{align*}
{
\gamma(t)
=
\hat\si(\hat X(t))^{-1}
\big[
\hat b(\hat X_t)-\hat b(\hat Y_t)
+\kappa(\hat X(t)-\hat Y(t))
\big],}
\qquad
\tt W(t)=W(t)+\int_0^t\gamma(s)\,\d s,
\end{align*}
and, for $t\geq0$, define
\begin{equation*}
R(t)
=
\exp\left(
-\int_0^t\langle\gamma(s),\,\d W(s)\rangle
-\ff12\int_0^t|\gamma(s)|^2\,\d s
\right).
\end{equation*}
We have the following results.
\begin{lem}\label{L3}
{Under the assumptions of Theorem \ref{T01},} if $\kappa>0$ is sufficiently large, then
\begin{align}\label{Rt}
\sup_{t\in[0,T]}\E(R(t)\log R(t))<\8,\ \ T>0.
\end{align}
Consequently, there exists a unique probability measure $\Q$ on $(\Omega,\F_\8)$ such that
\begin{align}\label{QP}
\ff{\d\Q|_{\F_t}}{\d\P|_{\F_t}}=R(t),\ \ t\geq0.
\end{align}
Moreover, $\tt W(t)$ is a $d$-dimensional Brownian motion under $\Q$.
\end{lem}
\begin{proof}
According to the proof of \cite[Lemma 3.2]{BWY19}, it suffices to prove \eqref{Rt}. For any $n\geq1$, define
\begin{equation*}
\tau_n=\inf\{t\geq0:\|\hat X_t\|_\tau+\|\hat Y_t\|_\tau\geq n\}.
\end{equation*}
By localization, the stopped exponential process
$(R(t\wedge\tau_n))_{t\geq0}$ is a martingale. For $T>0$, let
\begin{equation*}
\d\Q_{T,n}=R(T\wedge\tau_n)\d\P.
\end{equation*}
Then $(\tt W(t))_{t\in[0,T\wedge\tau_n]}$ is a Brownian motion under
$\Q_{T,n}$. Rewrite \eqref{dY} and \eqref{ttY} as
\begin{equation}\label{ttXY}
\begin{split}
&\d\hat X(t)
=
\{\hat b(\hat Y_t)-\kappa Z(t)\}\,\d t
+\hat\si(\hat X(t))\,\d\tt W(t),\\
&\d\hat Y(t)
=
\hat b(\hat Y_t)\,\d t
+\hat\si(\hat Y(t))\,\d\tt W(t),
\qquad t\leq T\wedge\tau_n,
\end{split}
\end{equation}
where $Z(t)=\hat X(t)-\hat Y(t)$. Hence,
\begin{equation}\label{Zeq}
\d Z(t)
=
-\kappa Z(t)\,\d t
+\{\hat\si(\hat X(t))-\hat\si(\hat Y(t))\}\,\d\tt W(t).
\end{equation}
Let $L>0$ satisfy
\begin{equation*}
\|\hat\si(x)-\hat\si(y)\|_{HS}\leq L|x-y|,
\qquad x,y\in\R^d.
\end{equation*}
For $j\geq2$, set
$
\rho_j:=\kappa-({j-1})L^2/{2}.
$
It\^o's formula gives
\begin{align*}
\d |Z(t)|^j
\leq&
-j\rho_j|Z(t)|^j\,\d t+j|Z(t)|^{j-2}
\left\langle Z(t),
\{\hat\si(\hat X(t))-\hat\si(\hat Y(t))\}\,\d\tt W(t)
\right\rangle .
\end{align*}
Choosing $\kappa$ such that $\rho_j>0$, applying It\^o's formula to
$\e^{j\rho_jt}|Z(t)|^j$ and stopping at $\tau_n$ yield
\begin{align}\label{Zpoint}
\E_{\Q_{T,n}}
[\1_{\{t\leq\tau_n\}}|Z(t)|^j]
\leq
\e^{-j\rho_jt}|Z(0)|^j.
\end{align}
Fix $\tau_1\in(0,\tau)$ and increase $\kappa$, if necessary, such that $\rho_j>\tau_1$. Put
$
a_j:=j(\rho_j-\tau_1)>0.
$
Applying It\^o's formula to $\e^{j\tau_1t}|Z(t)|^j$ gives
\begin{equation*}
\d\{\e^{j\tau_1t}|Z(t)|^j\}
\leq-a_j\e^{j\tau_1t}|Z(t)|^j\,\d t+\d M_j(t),
\end{equation*}
where
\begin{equation*}
M_j(t):=j\int_0^t\e^{j\tau_1s}|Z(s)|^{j-2}
\left\langle Z(s),
\{\hat\si(\hat X(s))-\hat\si(\hat Y(s))\}\,\d\tt W(s)
\right\rangle.
\end{equation*}
Consequently,
\begin{equation*}
a_j\E_{\Q_{T,n}}
\int_0^{t\wedge\tau_n}\e^{j\tau_1s}|Z(s)|^j\,\d s
\leq |Z(0)|^j.
\end{equation*}
Moreover, by the Burkholder--Davis--Gundy and Young inequalities, there exists a constant $c_5(j)>0$, independent of $T$ and $n$, such that
\begin{align*}
&\E_{\Q_{T,n}}
\left[
\sup_{s\in[0,t\wedge\tau_n]}\e^{j\tau_1s}|Z(s)|^j
\right]\\
&\quad\leq |Z(0)|^j
+c_5(j)\E_{\Q_{T,n}}
\left(
\int_0^{t\wedge\tau_n}\e^{2j\tau_1s}|Z(s)|^{2j}\,\d s
\right)^{1/2}\\
&\quad\leq |Z(0)|^j
+\frac12\E_{\Q_{T,n}}
\left[
\sup_{s\in[0,t\wedge\tau_n]}\e^{j\tau_1s}|Z(s)|^j
\right]
+c_5(j)\E_{\Q_{T,n}}
\int_0^{t\wedge\tau_n}\e^{j\tau_1s}|Z(s)|^j\,\d s.
\end{align*}
After enlarging $c_5(j)$, if necessary, we obtain
\begin{align}\label{Zsup}
\E_{\Q_{T,n}}
\left[
\sup_{s\in[0,t\wedge\tau_n]}
\e^{j\tau_1s}|Z(s)|^j
\right]
\leq
c_5(j)|Z(0)|^j.
\end{align}
Since $\tau_1<\tau$,
\begin{equation*}
\e^{\tau_1t}\|Z_t\|_\tau
\leq\|Z_0\|_\tau+
\sup_{s\in[0,t]}\e^{\tau_1s}|Z(s)|.
\end{equation*}
Hence, by \eqref{Zsup}, there exists a constant $c_6(j)>0$, independent of $T$ and $n$, such that
\begin{align}\label{Zseg}
\E_{\Q_{T,n}}
[\1_{\{t\leq\tau_n\}}\|Z_t\|_\tau^j]
\leq
c_6(j)\e^{-j\tau_1t}\|\xi-\eta\|_\tau^j.
\end{align}

By the definition of $\gamma(t)$ and \eqref{ttl}, there exists a constant $c_7>0$ such that
\begin{align}\label{gammabound}
|\gamma(t)|^2
&\leq
c_7\big\{
|\hat b(\hat X_t)-\hat b(\hat Y_t)|^2
+\kappa^2|Z(t)|^2
\big\}\notag\\
&\leq
c_7\big\{
(1+\kappa^2)\|Z_t\|_\tau^2
+\|\hat Y_t\|_\tau^{2\alpha}|Z(t)|^2
\big\}.
\end{align}
On the other hand, by \eqref{tt}, \eqref{btilde} and \eqref{nu2}, there exists a constant $c_8>0$ such that
\begin{equation*}
|\hat b(\zeta)|\leq c_8(1+\|\zeta\|_\tau),
\qquad
\|\hat\si\|_\8\leq c_8,
\qquad
\zeta\in\C_\tau.
\end{equation*}
By \eqref{EPR}, \eqref{ttXY},  the Burkholder--Davis--Gundy inequality and Gr\"onwall's inequality, for every $m\geq2$, there exists a constant $c_9(m)>0$, independent of $T$ and $n$,
such that
\begin{align}\label{Ymoment}
\E_{\Q_{T,n}}
\left[
\1_{\{t\leq\tau_n\}}\|\hat Y_t\|_\tau^m
\right]
\leq
c_9(m)\e^{c_9(m)t}(1+\|\eta\|_\tau^m),
\qquad t\in[0,T].
\end{align}

If $\alpha>0$, choose $j>2$ such that
$
m:={(2\alpha j)}/({j-2})\geq2.
$
By H\"older's inequality, \eqref{Ymoment} and \eqref{Zpoint}, there exists a constant $c_{10}>0$ such that
\begin{align*}
\E_{\Q_{T,n}}
[\1_{\{t\leq\tau_n\}}\|\hat Y_t\|_\tau^{2\alpha}|Z(t)|^2]
&\leq
\left(
\E_{\Q_{T,n}}
[\1_{\{t\leq\tau_n\}}\|\hat Y_t\|_\tau^m]
\right)^{\ff{j-2}{j}}
\left(
\E_{\Q_{T,n}}
[\1_{\{t\leq\tau_n\}}|Z(t)|^j]
\right)^{\ff2j}\\
&\leq
c_{10}(1+\|\eta\|_\tau^{2\alpha})
\e^{-c_\kappa t}|Z(0)|^2,
\end{align*}
where
$
c_\kappa
:=
2\rho_j-\ff{j-2}{j}c_9(m)>0
$
for sufficiently large $\kappa$. For $\alpha=0$, the same estimate holds with
$c_{10}=1$ and $c_\kappa=2\rho_2>0$ by \eqref{Zpoint} with $j=2$.

Since $|Z(0)|\leq\|\xi-\eta\|_\tau$, \eqref{Zseg}, \eqref{gammabound} and Tonelli's theorem imply that there exists a constant $c_{11}>0$, independent of $T$ and $n$, such that
\begin{equation}\label{gammaenergy}
\begin{aligned}
&\E_{\Q_{T,n}}
\int_0^{T\wedge\tau_n}|\gamma(s)|^2\,\d s
=
\int_0^T
\E_{\Q_{T,n}}
[\1_{\{s\leq\tau_n\}}|\gamma(s)|^2]\,\d s\\
&\quad\leq
c_{11}(1+\|\eta\|_\tau^{2\alpha})\|\xi-\eta\|_\tau^2.
\end{aligned}
\end{equation}
Consequently, for $t\in[0,T]$,
\begin{align*}
\E[R(t\wedge\tau_n)\log R(t\wedge\tau_n)]
&=
\ff12
\E_{\Q_{T,n}}
\int_0^{t\wedge\tau_n}|\gamma(s)|^2\,\d s\leq
\ff{c_{11}}{2}(1+\|\eta\|_\tau^{2\alpha})\|\xi-\eta\|_\tau^2.
\end{align*}
The uniform entropy estimate implies that
$\{R(t\wedge\tau_n)\}_{n\geq1}$ is uniformly integrable. Since
$
\lim_{n\to\infty}\tau_n=\infty
$
almost surely, it follows that $\E R(t)=1$ for every $t\geq0$, and hence
$R$ is a martingale. Fatou's lemma and the preceding bound prove
\eqref{Rt}. The consistent family of probability measures determined by
\eqref{QP} then defines $\Q$ on $\F_\infty$, and the last assertion follows
from Girsanov's theorem.
\end{proof}

\begin{rem}\label{dini}
{
\begin{enumerate}
\item[(1)] The Zvonkin transformation regularizes the Dini part $b^{(0)}$. The Lipschitz-type condition on $b^{(1)}$ yields the pointwise difference estimate in \eqref{ttl} for the transformed drift. Therefore,
$
\hat b(\hat X_t)-\hat b(\hat Y_t)
$
can be included in the change of probability measure via Girsanov's theorem, and this drift difference is cancelled in \eqref{Zeq}. This is different from \cite{BWY19}, where the drift is assumed to satisfy only a one-sided Lipschitz condition. Under that assumption, the square of the drift difference is not controlled by the coupling distance and hence cannot be included in the density in the same way; the drift difference remains in the equation for the coupling difference.

\item[(2)]
		The Dini condition on $b^{(0)}$ is used to obtain \eqref{nu1}, and hence
		$\hat\si$ is globally Lipschitz as stated in \eqref{ttl}. Together with the
		Lipschitz-type condition on $b^{(1)}$, this permits the coupling by change of
		measures used above, for which the drift difference is cancelled and
		\eqref{Zeq} holds. Under weaker Sobolev regularity, the same proof would
		require an exponential estimate for $G(\hat X)$ under $\Q$. Since $\hat X$
		is the controlled process, such an estimate does not follow directly from the
		usual Krylov estimates for the original equation. We therefore retain the
		Dini condition for the present method.
\end{enumerate}
}
\end{rem}

{\paragraph{$(a_3)$}In the next step, to deduce the asymptotic log-Harnack inequality from the asymptotic coupling $(\hat X_t,\hat Y_t)$ given in \eqref{ttXY}, we show that $\|\hat X_t-\hat Y_t\|_\tau$ decays exponentially fast as $t$ tends to $\8$ in the $L^p$-norm sense for all $p>0$.}

\begin{lem}\label{L4}
{Under the assumptions of Theorem \ref{T01}, for any $p>0$ and $\tau_0\in(0,\tau)$, there exist $\kappa,c>0$ such that}
\begin{align}\label{decay}
\E_\Q\|\hat X_t-\hat Y_t\|_\tau^p
\leq
c\e^{-p\tau_0t}\|\xi-\eta\|_\tau^p,\ \ t\geq0.
\end{align}
\end{lem}
\begin{proof}
{Fix $p>0$. Choose $j>\max\{p,2\}$ and $\tau_1\in(\tau_0,\tau)$.}
Increase $\kappa$, if necessary, such that $\rho_j>\tau_1$ and the estimates in Lemma \ref{L3} hold for this $j$. By Lemma \ref{L3}, \eqref{Zeq} holds under $\Q$ and $\tt W$ is a Brownian motion. Repeating the weighted It\^o--Burkholder--Davis--Gundy estimate leading to \eqref{Zsup}, followed by the segment estimate used in \eqref{Zseg}, gives
\begin{equation*}
\E_\Q\left[
\1_{\{t\leq\tau_n\}}\|Z_t\|_\tau^j
\right]
\leq
c_6(j)\e^{-j\tau_1t}\|\xi-\eta\|_\tau^j.
\end{equation*}
Since $\lim_{n\to\infty}\tau_n=\infty$ almost surely, the monotone convergence theorem yields
\begin{equation*}
\E_\Q\|Z_t\|_\tau^j
\leq
c_6(j)\e^{-j\tau_1t}\|\xi-\eta\|_\tau^j,
\qquad t\geq0.
\end{equation*}
Therefore, by Jensen's inequality,
\begin{align*}
\E_\Q\|Z_t\|_\tau^p
&\leq
\left(\E_\Q\|Z_t\|_\tau^j\right)^{p/j}\leq
c_6(j)^{p/j}\e^{-p\tau_1t}\|\xi-\eta\|_\tau^p
\leq
c_6(j)^{p/j}\e^{-p\tau_0t}\|\xi-\eta\|_\tau^p.
\end{align*}
Thus, \eqref{decay} holds with $c=c_6(j)^{p/j}$.
\end{proof}

{\paragraph{$(a_4)$}We now complete the proof of Theorem \ref{T01}.}
By Lemma \ref{L3} and the weak uniqueness of solutions to \eqref{dY}, we have
\begin{equation*}
\hat P_tf(\eta)=\E_\Q f(\hat Y_t),\qquad
t\geq0,\quad f\in\B_b(\C_\tau).
\end{equation*}
By Lemma \ref{L4}, there exists a constant $c_{12}>0$ such that, for any
$f\in\B_b^+(\C_\tau)$ with $\|\nn\log f\|_\8<\8$,
\begin{equation}\label{ttP1}
\begin{split}
\hat P_t\log f(\eta)
&=\E_\Q\log f(\hat Y_t)=\E_\Q\log f(\hat X_t)
+\E_\Q\{\log f(\hat Y_t)-\log f(\hat X_t)\}\\
&\leq
\E[R(t)\log f(\hat X_t)]
+\|\nn\log f\|_\8\E_\Q\|\hat X_t-\hat Y_t\|_\tau\\
&\leq
\E[R(t)\log R(t)]
+\log\hat P_tf(\xi)
+c_{12}\e^{-\tau_0t}
\|\nn\log f\|_\8\|\xi-\eta\|_\tau.
\end{split}
\end{equation}
Moreover, \eqref{gammaenergy} implies that there exists a constant $c_{13}>0$ such that
\begin{align*}
\E[R(t)\log R(t)]
&=
\ff12\E_\Q\int_0^t|\gamma(s)|^2\,\d s\leq
c_{13}(1+\|\eta\|_\tau^{2\alpha})\|\xi-\eta\|_\tau^2.
\end{align*}
Substituting this into \eqref{ttP1} gives the asymptotic log-Harnack inequality for $\hat P_t$. Together with \eqref{alp}, this completes the proof of Theorem \ref{T01}.

\section{  Proof of Theorem \ref{T02}}
\begin{proof}[Proof of Theorem \ref{T02}(1)]
Fix $T>0$, put $\vv:=\vv(\alpha)$, and let $\mu^1,\mu^2\in\scr P_{2+\vv,\e}^\alpha(\C_\tau)$. Consider
\begin{align}\label{Es}
\d X^i(t)
=
b(X^i_t,\mu^i_t)\,\d t+\si(X^i(t))\,\d W(t),
\quad t\in[0,T],\ i=1,2,
\end{align}
where $\mu^i_t=P_t^*\mu^i$, $i=1,2$.

\paragraph{$(\mathbf a)$}
According to \cite[Theorem 2.3]{WYZ}, these SDEs are well-posed. {Let $\delta_0\in\scr P(\C_\tau)$ be the Dirac measure at the zero path. Define $\bar b:\R^d\to\R^d$ and $\tt b:\C_\tau\times\scr P_2(\C_\tau)\to\R^d$ by
$
\bar b(x):=b^{(1)}(x^0,\delta_0), \ 
\tt b(\xi,\mu):=b^{(1)}(\xi,\mu)-\bar b(\xi(0)).
$
Then
\begin{align}\label{btildemu}
|\tt b(\xi,\mu)-\tt b(\eta,\nu)|
&\leq2K\|\xi-\eta\|_\tau+K_1\W_2(\mu,\nu),\notag\\
|\tt b(\xi,\mu)|
&\leq K(1+\|\xi\|_\tau^\alpha)+2K_1\|\mu\|_2.
\end{align}
Let $u$ solve \eqref{ee1} with this $\bar b$ and define $\Theta,\Theta^\#$ by \eqref{Theta}. Then
\begin{equation*}
\hat b(\xi,\mu)
=
\bar b(\Theta^{-1}(\xi(0)))+\lambda u(\Theta^{-1}(\xi(0)))
+\nn\Theta(\Theta^{-1}(\xi(0)))\tt b((\Theta^\#)^{-1}(\xi),\mu).
\end{equation*}
Moreover, there exists a constant $c_1>0$ such that
\begin{align}\label{ttlmu}
|\hat b(\xi,\mu)-\hat b(\eta,\nu)|
\leq&
c_1\|\xi-\eta\|_\tau+c_1K_1\W_2(\mu,\nu)
+c_1\{\|\eta\|_\tau^\alpha+K_1\|\nu\|_2\}
|\xi(0)-\eta(0)|.
\end{align}}
Let $Y^i(t)=\Theta(X^i(t))$, $i=1,2$, and denote $Z(t)=Y^1(t)-Y^2(t)$. Due to $(H)$, \eqref{nu1}, \eqref{nu2}, \eqref{ttlmu} and It\^o's formula, there exists a constant $c_2>0$ such that, for any $t\in[0,T]$,
\begin{equation*}
\begin{aligned}
\d\e^{-A(t)}
\left(\e^{\tau t}|Z(t)|\right)^{2+\vv}
\leq&
c_2\e^{-A(t)}
\left\{
\left(\e^{\tau t}\|Z_t\|_\tau\right)^{2+\vv}
+\left(\e^{\tau t}\W_2(\mu_t^1,\mu_t^2)\right)^{2+\vv}
\right\}\,\d t+\d\tt M(t),
\end{aligned}
\end{equation*}
where
\begin{equation}\label{am}
\begin{split}
A(t)
&=
c_2\int_0^t
\{\|X_s^2\|_\tau^\alpha+K_1\|\mu_s^2\|_2\}\,\d s,\\
\d\<\tt M\>(t)
&\leq
c_2\e^{-2A(t)}
\left(\e^{\tau t}|Z(t)|\right)^{2(2+\vv)}\,\d t.
\end{split}
\end{equation}
Thus, {since $A(t)$ is increasing}, \eqref{EPR}, the Burkholder--Davis--Gundy and Young inequalities yield a constant $c_3>0$ such that
\begin{equation}\label{u}
\begin{split}
U_t
:=&
\E\left[
\sup_{s\in[0,t]}
\e^{-A(s)}
\left(\e^{\tau s}\|Z_s\|_\tau\right)^{2+\vv}
\bigg|\F_0
\right]\\
\leq&
c_3\e^{c_3t}
\left\{
\|Z_0\|_\tau^{2+\vv}
+\int_0^t
\e^{(2+\vv)\tau s}
\W_2(\mu_s^1,\mu_s^2)^{2+\vv}\,\d s
\right\}.
\end{split}
\end{equation}
By \cite[Corollary 4.5]{WYZ} and \cite[Theorem 2.3]{WYZ}, for any $\beta,\dd>0$ there exists a constant $c_4(T,\beta,\dd)>0$ such that
\begin{align}\label{eam}
\E[\e^{\beta A(t)}|\F_0]
\leq
c_4(T,\beta,\dd)
\e^{\dd\|X_0^2\|_\tau^{2\alpha}},
\qquad t\in[0,T].
\end{align}
When $\alpha>0$ and $\vv>0$, choose $\dd\in(0,1)$ such that $\mu^2(\e^{\dd\|\cdot\|_\tau^{2\alpha}})<\8$. Then \eqref{u}, \eqref{eam} and H\"older's inequality give a constant $c_5(T)>0$ such that
\begin{equation*}
\begin{split}
&\E\left[
\sup_{s\in[0,t]}
\e^{2\tau s}\|Z_s\|_\tau^2
\bigg|\F_0
\right]\leq
c_5(T)
\e^{\ff{\vv\dd}{2+\vv}\|X_0^2\|_\tau^{2\alpha}}
\left\{
\|Z_0\|_\tau^2
+
\left(
\int_0^t
\e^{(2+\vv)\tau s}
\W_2(\mu_s^1,\mu_s^2)^{2+\vv}\,\d s
\right)^{\ff2{2+\vv}}
\right\}.
\end{split}
\end{equation*}

\paragraph{$(\mathbf b)$}
By \eqref{Wfull}, take $\F_0$-measurable random variables $X_0^1,X_0^2$ with distributions $\mu^1,\mu^2$ such that
\begin{equation}\label{M1}
\E\|X_0^1-X_0^2\|_\tau^{2+\vv}
=
\W_{2+\vv}(\mu^1,\mu^2)^{2+\vv}.
\end{equation}
Let $X^i$ be the corresponding solutions and set
\begin{equation*}
Z(t):=\Theta(X^1(t))-\Theta(X^2(t)),\qquad t\in[0,T].
\end{equation*}
When $\alpha>0$ and $\vv>0$, taking expectations in the preceding conditional estimate and using H\"older's inequality, \eqref{M1} and \eqref{nu2}, there exists a constant $c_6(T)>0$ such that
\begin{align*}
\W_2(\mu_t^1,\mu_t^2)^{2+\vv}
&\leq
2^{2+\vv}\left(\E\|Z_t\|_\tau^2\right)^{\ff{2+\vv}{2}}\\
&\leq
c_6(T)
\left(\mu^2(\e^{\dd\|\cdot\|_\tau^{2\alpha}})\right)^{\ff\vv2}
\left\{
\W_{2+\vv}(\mu^1,\mu^2)^{2+\vv}
+\int_0^t\W_2(\mu_s^1,\mu_s^2)^{2+\vv}\,\d s
\right\}.
\end{align*}
By Gr\"onwall's lemma, there exists a constant $c_7(T)>0$ such that
\begin{align}\label{wfinite}
\W_2(\mu_t^1,\mu_t^2)
\leq c_7(T)\W_{2+\vv}(\mu^1,\mu^2),
\qquad t\in[0,T].
\end{align}
If $\alpha=0$, taking $\vv=0$ in \eqref{u} gives, for some constant $c_8(T)>0$,
\begin{align*}
\W_2(\mu_t^1,\mu_t^2)^2
&\leq
4\E\|Z_t\|_\tau^2\\
&\leq
c_8(T)
\left\{
\W_2(\mu^1,\mu^2)^2
+\int_0^t\W_2(\mu_s^1,\mu_s^2)^2\,\d s
\right\},
\end{align*}
and \eqref{wfinite} follows again by Gr\"onwall's lemma, after enlarging $c_7(T)$ if necessary.
The constant $c_7(T)$ in \eqref{wfinite} may be chosen nondecreasing and locally bounded in $T$. Defining $c_0(t):=c_7(t)$, {the arbitrariness of $T>0$ proves \eqref{wt}.}
\end{proof}

Now, let $\tau_0\in(0,\tau)$. It remains to verify \eqref{log1}, which implies \eqref{lip} through repeating the proof of \cite[Theorem 2.1(1)]{BWY19}.

\beg{proof}[Proof of \eqref{log1}]
The case $t=0$ is immediate. Fix $t>0$ and
$\mu,\nu\in\scr P_{2+\vv,\e}^\alpha(\C_\tau)$, where
$\vv=\vv(\alpha)$. By \eqref{Wfull}, take $\F_0$-measurable random variables
$X_0^\mu,X_0^\nu$ with distributions $\mu,\nu$ such that
\begin{equation}\label{M3}
\E\|X_0^\mu-X_0^\nu\|_\tau^{2+\vv}
=
\W_{2+\vv}(\mu,\nu)^{2+\vv}.
\end{equation}
Let $X^\mu$ and $X^\nu$ be the solutions to \eqref{E-1} with initial values
$X_0^\mu$ and $X_0^\nu$, respectively, and set
$\mu_s=P_s^*\mu$ and $\nu_s=P_s^*\nu$ for $s\in[0,t]$. Let
\begin{align*}
&\bar\zeta_s
:=
\{\si^*(\si\si^*)^{-1}\}(X^\mu(s))
[b(X^\mu_s,\mu_s)-b(X^\mu_s,\nu_s)],\\
&\bar R_s
:=
\exp\left[
-\int_0^s\langle\bar\zeta_u,\,\d W(u)\rangle
-\ff12\int_0^s|\bar\zeta_u|^2\,\d u
\right],
\qquad s\in[0,t].
\end{align*}
By $(H)$ and \eqref{wt}, there exists a constant $c_2>0$ such that
\begin{align}\label{bzt}
|\bar\zeta_s|
\leq
c_2K_1\W_2(\mu_s,\nu_s)
\leq
c_2K_1c_0(s)\W_{2+\vv}(\mu,\nu).
\end{align}
Thus, $\bar R_s$ is a martingale and
\begin{equation*}
\bar W(s)=W(s)+\int_0^s\bar\zeta_r\,\d r
\end{equation*}
is a Brownian motion under $\bar\P:=\bar R_t\P$. Then
\begin{align}\label{xmu}
\d X^\mu(s)
=
b(X^\mu_s,\nu_s)\,\d s+\si(X^\mu(s))\,\d\bar W(s),
\qquad s\in[0,t].
\end{align}
Let $Y^\mu(s):=\Theta(X^\mu(s))$. By It\^o's formula,
\begin{align}\label{ymu}
\d Y^\mu(s)
=
\hat b(Y^\mu_s,\nu_s)\,\d s
+\hat\si(Y^\mu(s))\,\d\bar W(s),
\qquad s\in[0,t],
\end{align}
where $\hat b$ is given in the proof of Theorem \ref{T02}(1).

Since $\E(\bar R_t|\F_0)=1$, the distribution of
$(X_0^\mu,X_0^\nu)$ is unchanged under $\bar\P$. Under $\bar\P$, let
$\hat Y_0=\Theta^\#(X_0^\nu)$ and consider the unique strong solution of
\begin{equation}\label{coupling2}
\begin{aligned}
\d\hat Y(s)
=&
\Big\{
\hat b(\hat Y_s,\nu_s)
+\hat\si(\hat Y(s))\hat\si(Y^\mu(s))^{-1}\\
&\quad\times
[
\hat b(Y^\mu_s,\nu_s)-\hat b(\hat Y_s,\nu_s)
+\kappa(Y^\mu(s)-\hat Y(s))
]
\Big\}\,\d s+\hat\si(\hat Y(s))\,\d\bar W(s).
\end{aligned}
\end{equation}
For $s\in[0,t]$, let
\begin{equation*}
\tt\zeta_s
=
\hat\si(Y^\mu(s))^{-1}
[
\hat b(Y^\mu_s,\nu_s)-\hat b(\hat Y_s,\nu_s)
+\kappa(Y^\mu(s)-\hat Y(s))
]
\end{equation*}
and
\begin{equation*}
\tt R_s
=
\exp\left[
-\int_0^s\langle\tt\zeta_u,\,\d\bar W(u)\rangle
-\ff12\int_0^s|\tt\zeta_u|^2\,\d u
\right].
\end{equation*}
As in Lemma \ref{L3}, define
\begin{equation*}
\sigma_n
:=
\inf\{s\geq0:\|Y^\mu_s\|_\tau+\|\hat Y_s\|_\tau\geq n\},
\qquad n\geq1,
\end{equation*}
and let
\begin{equation*}
\d\Q_n
=
\tt R_{t\wedge\sigma_n}\,\d\bar\P
\end{equation*}
on $\F_t$. Then
\begin{equation*}
\tt W(s)=\bar W(s)+\int_0^s\tt\zeta_r\,\d r
\end{equation*}
is a Brownian motion up to $t\wedge\sigma_n$ under $\Q_n$, and
\begin{align*}
\d Y^\mu(s)
&=
\{\hat b(\hat Y_s,\nu_s)-\kappa Z(s)\}\,\d s
+\hat\si(Y^\mu(s))\,\d\tt W(s),\\
\d\hat Y(s)
&=
\hat b(\hat Y_s,\nu_s)\,\d s
+\hat\si(\hat Y(s))\,\d\tt W(s),
\end{align*}
where $Z(s)=Y^\mu(s)-\hat Y(s)$. Thus, up to $\sigma_n$, the difference
process satisfies the same equation as \eqref{Zeq}. Repeating the proof of
\eqref{Zpoint}--\eqref{Zseg} conditionally on $\F_0$, and using \eqref{M3}
and \eqref{nu2}, we may choose $\kappa$ sufficiently large such that, for
some constant $c_3>0$ independent of $n$,
\begin{equation*}
\E_{\Q_n}
\left[
\1_{\{s\leq\sigma_n\}}\|Y^\mu_s-\hat Y_s\|_\tau^2
\right]
\leq
c_3\e^{-2\tau_0s}\W_{2+\vv}(\mu,\nu)^2,
\qquad s\in[0,t].
\end{equation*}

Moreover, by \eqref{ttlmu}, there exists a constant $c_4>0$ such that
\begin{equation*}
|\tt\zeta_s|^2
\leq
c_4\left\{
\|Z_s\|_\tau^2+
\big(\|\hat Y_s\|_\tau^{2\alpha}
+K_1^2\|\nu_s\|_2^2\big)|Z(s)|^2
\right\}.
\end{equation*}
By \cite[Theorem 2.3]{WYZ} and repeating the argument leading to
\eqref{Ymoment}--\eqref{gammaenergy}, with the additional term involving
$\|\nu_s\|_2$ controlled by \eqref{Zpoint}, we may increase $\kappa$ if
necessary such that, for some constant $c_5>0$ independent of $t$ and $n$,
\begin{equation*}
\sup_{n\geq1}
\E_{\Q_n}
\int_0^{t\wedge\sigma_n}|\tt\zeta_s|^2\,\d s
\leq
c_5\W_{2+\vv}(\mu,\nu)^2.
\end{equation*}
Consequently,
\begin{align*}
\E_{\bar\P}
[\tt R_{t\wedge\sigma_n}\log\tt R_{t\wedge\sigma_n}]
&=
\ff12\E_{\Q_n}
\int_0^{t\wedge\sigma_n}|\tt\zeta_s|^2\,\d s\leq
\ff{c_5}{2}\W_{2+\vv}(\mu,\nu)^2.
\end{align*}
Thus, the stopped densities are uniformly integrable. Since
$
\lim_{n\to\8}\sigma_n=\8
$
almost surely, the stopping can be removed as in Lemma \ref{L3}. Hence,
$\tt R_s$ is a martingale, $\Q:=\tt R_t\bar\P$ is well defined, and
$\tt W$ is a Brownian motion under $\Q$. Letting $n$ tend to infinity in
the preceding stopped estimates yields
\begin{align}\label{M4}
\E_\Q\|Y^\mu_s-\hat Y_s\|_\tau^2
\leq
c_3\e^{-2\tau_0s}\W_{2+\vv}(\mu,\nu)^2,
\qquad s\in[0,t],
\end{align}
and
\begin{align}\label{tzeta}
\E_\Q\int_0^t|\tt\zeta_s|^2\,\d s
\leq
c_5\W_{2+\vv}(\mu,\nu)^2,
\qquad t\geq0.
\end{align}

Since
\begin{equation*}
\E(\bar R_t|\F_0)=\E_{\bar\P}(\tt R_t|\F_0)=1,
\end{equation*}
the initial joint distribution is unchanged under $\Q$. By weak uniqueness
of the transformed equation,
\begin{equation*}
\L_{\hat Y_s|\Q}
=
\L_{\Theta^\#(X_s^\nu)}.
\end{equation*}
By \eqref{nu2}, $(\Theta^\#)^{-1}$ is $2$-Lipschitz. Hence, applying
Young's inequality as in \cite[Lemma 2.4]{ATW} and using \eqref{M4}, there
exists a constant $c_6>0$ such that
\begin{equation}\label{log3}
\begin{split}
P_t\log f(\nu)
&=
\E_\Q\log(f\circ(\Theta^\#)^{-1})(\hat Y_t)\\
&\leq
\E_\Q\log(f\circ(\Theta^\#)^{-1})(Y^\mu_t)
+2\|\nn\log f\|_\8\E_\Q\|Y^\mu_t-\hat Y_t\|_\tau\\
&\leq
\E[R_t\log R_t]
+\log P_tf(\mu)
+c_6\e^{-\tau_0t}\|\nn\log f\|_\8
\W_{2+\vv}(\mu,\nu),
\end{split}
\end{equation}
where $R_t:=\bar R_t\tt R_t$ is the density of $\Q$ with respect to $\P$
on $\F_t$. Since
\begin{equation*}
R_t
=
\exp\left[
-\int_0^t\langle\bar\zeta_s+\tt\zeta_s,\,\d W(s)\rangle
-\ff12\int_0^t|\bar\zeta_s+\tt\zeta_s|^2\,\d s
\right],
\end{equation*}
\eqref{bzt}, \eqref{tzeta} and Young's inequality yield a constant
$c_7>0$ such that
\begin{align*}
\E[R_t\log R_t]
&=
\ff12\E_\Q\int_0^t|\bar\zeta_s+\tt\zeta_s|^2\,\d s\\
&\leq
c_7\left\{
1+K_1^2\int_0^t c_0(s)^2\,\d s
\right\}
\W_{2+\vv}(\mu,\nu)^2\leq
c_1(t)\W_{2+\vv}(\mu,\nu)^2.
\end{align*}
Putting this back into \eqref{log3} yields \eqref{log1}.
\end{proof}
\section*{Statements and Declarations}
\noindent\textbf{Acknowledgement:} The author would like to thank the referees for their corrections and valuable comments.

\noindent\textbf{Funding:} This work is supported in part by  the National Key R\&D Program of China (No. 2022YFA1006000), NNSFC (12531007).

\noindent\textbf{Conflict of interest:} The author states that there is no conflict of interest.

\noindent\textbf{Data availability:} {Data sharing is not applicable to this article, as no datasets were generated or analysed during the current study.}
\bibliographystyle{plain}
\bibliography{reference}
\end{document}